\RequirePackage{fix-cm}
\documentclass[smallextended]{svjour3}       
\smartqed  
\usepackage{graphicx}
\usepackage[utf8]{inputenc}
\usepackage[T1]{fontenc}
\usepackage{lmodern}
\usepackage[a4paper]{geometry}
\usepackage[english]{babel}
\usepackage{amsmath,amssymb}
\usepackage{version} 
\usepackage{color}

\newcommand{\Miln}{\operatorname{Miln}}

\newcommand{\Hyp}{\operatorname{Hyp}}

\newcommand{\dpar}{\partial}
\newcommand{\tD}{\tilde D}
\newcommand{\tH}{\tilde H}
\newcommand{\iso}{\overset{\sim}{\longrightarrow}}
\newcommand{\isom}{\overset{\sim}{=}}

\newcommand{\Spec}{\operatorname{Spec}}

\newcommand{\fgl}{\mathfrak{gl}}
\newcommand{\fsl}{\mathfrak{sl}}
\newcommand{\fsu}{\mathfrak{su}}

\newcommand{\Exp}{\operatorname{Exp}}

\newcommand{\Diag}{\operatorname{Diag}}

\newcommand{\be}{\begin{equation}}
\newcommand{\ee}{\end{equation}}

\newcommand{\lra}{\longrightarrow}

\newcommand{\BC}{\mathbb{C}}
\newcommand{\BR}{\mathbb{R}}

\newcommand{\fg}{\mathfrak g}

\newcommand{\fk}{\mathfrak k}

\begin{document}

\title{Vanishing cycles and Cartan eigenvectors}

\author{Laura Brillon \and Revaz Ramazashvili\footnote{R. R. contributed  only to writing the Section 6.} 
\and
	Vadim Schechtman \and Alexander Varchenko
}

\institute{L. B., V. S. : \at
              Institut de Mathématiques de Toulouse, 118 route de Narbonne, 31062 Toulouse, France \\
              \email{laura.brillon@math.univ-toulouse.fr , vadim.schechtman@math.univ-toulouse.fr}
              \and R.R. : \at Laboratoire de Physique Th\'eorique, IRSAMC, 
Universit\'e de Toulouse, CNRS, UPS, 31062 Toulouse, France \\
\email{revaz.ramazashvili@irsamc.ups-tlse.fr}       
           \and
           A. V. : \at
               Mathematics Department, University of North Carolina, Chapel Hill, NC, USA \\
               \email{anv@email.unc.edu} 
}

\date{Received: date / Accepted: date}

\maketitle


\centerline{\today}

\bigskip\bigskip

\begin{abstract}
Using the ideas coming from the singularity theory, we study the eigenvectors 
of the Cartan matrices of finite root systems, and of  $q$-deformations 
of these matrices. 
	


\end{abstract}

\section{Introduction}
\label{intro}

Let $A(R)$ be the Cartan matrix of  a finite root 
system $R$. The coordinates of its eigenvectors have an important 
meaning in the physics of integrable systems; we will say more on this below. 

The aim of this note is  to study these numbers and 
their $q$-deformations, using some results coming from the singularity theory.

\bigskip

We discuss three ideas:

(a) Cartan/Coxeter correspondence;

(b) Sebastiani - Thom product;

(c) Givental's $q$-deformations.

Let us explain what we are talking about.

\bigskip

Let us suppose that $R$ is simply laced, i.e. of type $A, D$, or $E$. 
These root systems are in one-to-one correspondence with (classes of) simple singularities 

$f: \mathbb{C}^N \to \mathbb{C}$, cf. \cite{AGV}. Under this correspondence, the root 
lattice $Q(R)$ is identified with the lattice of vanishing cycles, and 
the Cartan matrix $A(R)$ is the intersection matrix with respect to  
a {\it distinguished base}. The action of the Weyl group on 
$Q(R)$ is realized by Gauss - Manin monodromies - this is the Picard - Lefschetz theory (for some details see \S\ref{sec2} below). 

\bigskip

Remarkably, this  geometric picture provides a finer structure: namely, the symmetric matrix $A = A(R)$ comes equipped with a decomposition

\begin{equation}
A = L + L^t
\end{equation}

where $L$ is  a nondegenerate triangular ''Seifert form'', or ''variation matrix''. The matrix

\begin{equation}
C = - L^{-1}L^t
\end{equation}

represents a Coxeter element of $R$; geometrically it is 
the operator of ''classical monodromy''.

\bigskip 

We call the relation (1) - (2) between the Cartan matrix and the Coxeter element the {\it Cartan/Coxeter correspondence}. It works 
more generally for non-symmetric $A$ (in this case  (1) should be replaced by
 
\begin{equation}\label{non-sym}
A = L + U
\end{equation}

where $L$ is lower triangular and $U$ is upper triangular), and is due to Coxeter, 
cf. \cite{Co}, no. 1, p. 767, see \S 3 below.  

\bigskip

In a particular case (corresponding to a bipartition of the Dynkin 
graph) this relation is equivalent to an observation by R.Steinberg, 
cf. \cite{Stein}, cf. \S \ref{sub34} below.

\bigskip
 
This corresppondence allows one to relate the eigenvectors of $A$ and $C$, cf. Theorem \ref{vp}.

\bigskip

A decomposition (1) will be called {\it a polarization} of the 
Cartan matrix $A$. In \ref{sub33} below we introduce an operation of 
{\it Sebastiani - Thom}, or {\it joint} product $A*B$ of 
Cartan matrices (or of polarized lattices) $A$ and $B$. The root lattice of $A*B$ is the tensor product of the root lattice of $A$ 
and the root lattice of $B$.
With respect to this operation the Coxeter eigenvectors factorize very simply.

\bigskip

For example, the lattices $E_6$ and $E_8$ decompose into three ''quarks'': 

\begin{equation}\label{e6-dec}
E_6 = A_3*A_2*A_1
\end{equation}

\begin{equation}\label{e8-dec}
E_8 = A_4*A_2*A_1
\end{equation}

These decompositions are the main message from the singularity theory, 
and we discuss them in detail in this note. 

We use (\ref{e6-dec}), (\ref{e8-dec}), and the Cartan/Coxeter correspondence to get  expressions for all Cartan eigenvectors of $E_6$ and $E_8$; this is the first main result of this note, 
see \ref{sub491}, \ref{sub-E6} below.

(An elegant expression for all the Cartan eigenvectors of all finite root systems was given  by P.Dorey, cf. \cite{D} (a), Table 2 on p. 659.)   

\bigskip

In the paper \cite{Giv}, A. Givental has proposed a $q$-twisted version 
of the Picard - Lefschetz theory, which gave rise to a $q$-deformation of $A$,

\begin{equation}
A(q) = L + qL^t.
\end{equation}

Again, as Givental remarked, the decomposition (\ref{non-sym}) 
allows us to drop the assumption of symmetry in the definition above. 
In the last section, \S \ref{sec5}, we calculate the eigenvalues and eigenvectors 
of $A(q)$ in terms of the eigenvalues and eigenvectors of $A$. 
This is the second main result of this note. 

It turns out that if $\lambda$ is an eigenvalue of $A$ then 

\begin{equation}
\lambda(q) = 1 + (\lambda - 2)\sqrt{q} + q
\end{equation}

will be an eigenvalue of $A(q)$.  The coordinates of the corresponding 
eigenvector $v(q)$ are obtained from the coordinates of $v = v(1)$ 
by multiplication by appropriate powers of $q$; this is related to the fact that the Dynkin graph of $A$ is a tree, cf.  \ref{sub-tree}. For an example of $E_8$, see (\ref{eq1}).

\bigskip

In physics, the coordinates of the Perron - Frobenius Cartan eigenvectors 
appear as  particle masses  in affine Toda field theories, cf. \cite{D,F} 
and the Section \ref{toda} below. 

In a pioneering paper \cite{Z}, A. B. Zamolodchikov has discovered 
an octuplet of particles of $E_8$ symmetry in 
the two-dimensional critical Ising model in a magnetic field, 
and calculated their masses, cf. the Sections 
\ref{sub492} and \ref{sub64}.


The Appendix outlines some of the results of a neutron scattering 
experiment \cite{Coldea}, where the two lowest-mass $E_8$ particles of 
the Zamolodchikov's theory may have been observed.  

\section{Recollections from the singularity theory}
\label{sec2}

Here we recall some classical constructions and statements, cf. \cite{AGV}. 


\subsection{Lattice of vanishing cycles}
\label{sub21}

 Let $f:\ (\mathbb{C}^N,0) \to (\mathbb{C},0)$ be the germ of  a holomorphic function with an isolated 
 critical point at $0$, with $f(0) = 0$. We will be interested only in polynomial 
 functions (from the list below, cf. \S \ref{sub24}), so $f\in \mathbb{C}[x_1,\ldots, x_N]$. The {\it Milnor ring} 
 of $f$ is defined by 
 $$
 \Miln(f,0) = \mathbb{C}[[x_1,\ldots, x_N]]/(\dpar_1f,\ldots, \dpar_Nf)
 $$
 where $\dpar_i := \dpar/\dpar x_i$; it is a finite-dimensional 
 commutative $\mathbb{C}$-algebra. (In fact, it is a Frobenius, or, equivalently, a Gorenstein algebra.) The number
 $$
 \mu := \dim_\mathbb{C} \Miln(f,0)
 $$
 is called  the multiplicity or Milnor number of $(f,0)$. 
 
 A {\it Milnor fiber} is 
 $$
 V_z = f^{-1}(z)\cap \bar B_\rho
 $$
 where 
 $$
 \bar B_\rho = \{(x_1,\ldots, x_N)|\ \sum |x_i|^2 \leq \rho\}
 $$
 for  $1\gg \rho\gg |z|>0$.
 
 For $z$ belonging to a small disc $D_\epsilon = \{z\in \mathbb{C}|\ |z| < \epsilon\}$, the space $V_z$  is a complex manifold with boundary, homotopically 
 equivalent to a bouquet $\vee S^{N-1}$ 
 of $\mu$ spheres, \cite{M}. 
 
 The family  of free abelian groups 
 
 \begin{equation}
  Q(f;z) := \tH_{N-1}(V_z;\mathbb{Z})\isom \mathbb{Z}^\mu,\ z\in \overset{\bullet}D_\epsilon := 
  D_\epsilon \setminus \{0\}, 
 \end{equation}

 ($\tH$ means that we take  the reduced homology for $N = 1$), 
 carries a flat Gauss - Manin conection. 
 
 Take $t\in\mathbb{R}_{>0}\cap \overset{\bullet}D_\epsilon$; the lattice 
 $Q(f;t)$ does not depend, up to a canonical isomorphism, 
 on the choice of $t$. Let us call this lattice $Q(f)$. The linear operator 
 
 \begin{equation}
 T(f): Q(f) \iso Q(f)
 \end{equation}

 induced by the path $p(\theta) = e^{i\theta}t,\ 0\leq \theta\leq 2\pi$, is called the classical monodromy of the germ $(f,0)$. 
 
 In all the examples below $T(f)$ has finite order $h$.  The  
 eigenvalues of $T(f)$ have the form $e^{2\pi i k/h},\ k\in \mathbb{Z}$.  The set of suitably chosen $k$'s for each eigenvalue are called  
 the {\it spectrum} of our singularity.  

\subsection{Morse deformations}
\label{sub22}

 The $\mathbb{C}$-vector space $\Miln(f,0)$ may be identified with the 
 tangent space to the base $B$ of the miniversal defomation 
 of $f$. For 
 $$
 \lambda \in B^0 = B\setminus \Delta
 $$ 
 where $\Delta\subset B$ is an analytic subset of codimension $1$, 
 the corresponding function $f_\lambda: \mathbb{C}^N\to \mathbb{C}$ has $\mu$ nondegenerate Morse critical points with distinct critical values, 
 and the algebra $\Miln(f_\lambda)$ is semisimple, 
 isomorphic to $\mathbb{C}^\mu$. 
 
 Let $0\in B$ denote the point corresponding to $f$ itself, so that 
 $f = f_0$, and pick $t\in \mathbb{R}_{>0}\cap \overset{\bullet}D_\epsilon$ as in \S \ref{sub21}. 
 
 Afterwards  
 pick $\lambda\in B^0$ close to $0$ in such a way that the critical values $z_1, \ldots z_\mu$ of $f_\lambda$ have absolute values $\ll t$.  
 
 As in \S \ref{sub21},  
 for each 
 $$
 z\in \tD_\epsilon := D_\epsilon\setminus\{z_1, \ldots z_\mu\}
 $$
 the Milnor fiber $ V_z$ has the homotopy type of  
 a bouquet $\vee S^{N-1}$ 
 of $\mu$ spheres, and we will be interested in the middle 
 homology 
 $$
 Q(f_\lambda;z) = \tH_{N-1}(V_z;\mathbb{Z})\isom \mathbb{Z}^\mu
 $$
 
 The lattices $Q(f_\lambda;z)$ carry a natural bilinear  product induced by the  cup product in the homology which is symmetric (resp. skew-symmetric) when $N$ is 
 odd (resp. even).     
 
 The collection of these lattices, 
 when $z\in \tD_\epsilon$ varies, carries a flat Gauss - Manin connection. 
 
 Consider an ''octopus''  
 $$
 Oct(t)\subset\mathbb{C}
 $$ 
 with the head at $t$: a collection of non-intersecting 
 paths $p_i$ (''tentacles'') connecting $t$ with $z_i$ and not meeting the critical values $z_j$ otherwise. It gives rise 
 to a base 
 $$
 \{b_1,\ldots, b_\mu\}\subset Q(f_\lambda) := Q(f_\lambda;t)
 $$ 
 (called ''distinguished'') where $b_i$ is the cycle vanishing  
 when being transferred from $t$ to $z_i$ along the tentacle $p_i$,  
 cf. \cite{Gab}, \cite{AGV}. 
 
 The Picard - Lefschetz formula  describes  the action of the fundamental group $\pi_1(\tD_\epsilon;t)$ on 
 $Q(f_\lambda)$ with respect to this basis. Namely, 
 consider a loop $\gamma_i$ which turns around $z_i$ along the tentacle $p_i$,  
 then the corresponding transformation of $Q(f_\lambda)$ is the reflection 
 (or transvection) $s_i := s_{b_i}$, cf. \cite{Lef}, Th\'eor\`eme fondamental, Ch. II, p. 23. 
 
 The loops $\gamma_i$ generate the fundamental group 
 $\pi_1(\tD_\epsilon)$. Let 
 $$
 \rho:\ \pi_1(\tD_\epsilon;t)\to GL(Q(f_\lambda))
 $$
 denote the monodromy representation. The image of $\rho$, denoted by  $G(f_\lambda)$ and called the {\it monodromy group of $f_\lambda$}, 
 lies inside the subgroup 
 \newline $O(Q(f_\lambda))\subset GL(Q(f_\lambda))$ 
 of linear transformations respecting the above mentioned 
 bilinear form on $Q(f_\lambda)$.    
 
 The subgroup $G(f_\lambda)$ is generated by $s_i, 1\leq i \leq \mu$. 
 
 As in \S \ref{sub21}, we have the monodromy operator 
 $$
 T(f_\lambda)\in G(f_\lambda),
 $$
 the image by $\rho$ of the path $p\subset \tD_\epsilon$ starting at $t$ and going around all points $z_1, \ldots, z_\mu$.
 
 This operator 
 $T(f_\lambda)$ is now a product of $\mu$ simple reflections
 $$
 T(f_\lambda) = s_1s_2\ldots s_\mu,
 $$
 - this is because the only critical value $0$ of $f$ became 
 $\mu$ critical values $z_1, \ldots, z_\mu$ of $f_\lambda$.  
 
 One can identify the relative (reduced) homology
 $\tH_{N-1}(V_t, \dpar V_t;\mathbb{Z})$ with the dual group 
 $\tH_{N-1}(V_t;\mathbb{Z})^*$,  and one defines a map 
 $$
 \text{var}: \tH_{N-1}(V_t, \dpar V_t;\mathbb{Z}) \to \tH_{N-1}(V_t;\mathbb{Z}),
 $$
 called a {\it variation operator}, which translates to a map
 $$
 L: Q(f_\lambda)^* \iso Q(f_\lambda)
 $$
 (''Seifert  form'') 
 such that the matrix $A(f_\lambda)$ of the bilinear form 
 in the distinguished basis is 
 $$
 A(f_\lambda) = L + (-1)^{N-1}L^t,
 $$
 and 
 $$
 T(f_\lambda) = (-1)^{N-1}LL^{-t}.
 $$
 

 A choice of a path $q$ in $B$ connecting $0$ with $\lambda$, 
 enables one to identify $Q(f)$ with $Q(f_\lambda)$, and 
 $T(f)$ will be identified with $T(f_\lambda)$. 
 
 The image $G(f)$ of the monodromy group $G(f_\lambda)$ in 
 $GL(Q(f)) \isom GL(Q(f_\lambda))$ is  called the monodromy group 
 of $f$; it does not depend on a choice of a path $q$.
 
 \subsection{Sebastiani - Thom factorization}
\label{sub23}

If $g\in \mathbb{C}[y_1,\ldots, y_M]$ is another function, the sum, or 
{\bf join} of two singularities 
$f\oplus g:\ \mathbb{C}^{N+M}\to \mathbb{C}$ is defined by 
$$
(f\oplus g)(x,y) = f(x) + g(y)
$$
Obviously we can identify 
$$
\Miln(f\oplus g)\isom \Miln(f)\otimes \Miln(g)
$$
Note that the function $g(y) = y^2$ is a unit for this operation.

It follows that the singularities 
$f(x_1,\ldots, x_N)$ and 
$$
f(x_1,\ldots, x_N) + x_{M+1}^2 + \ldots + x^2_{N+M}
$$
are ''almost the same''. In order to have good signs (and for other purposes) it is 
convenient to add some squares to a given $f$ to get 
$N\equiv 3\mod(4)$.  

The fundamental Sebastiani - Thom theorem, \cite{ST}, says that there exists a natural isomorphism of lattices 
$$
Q(f\oplus g) \isom Q(f)\otimes_\mathbb{Z} Q(g),
$$
and under this identification 
the full monodromy decomposes as 
$$
T_{f\oplus g} = T_f\otimes T_g
$$
Thus, if 
$$
\Spec(T_f) = \{e^{\mu_p\cdot 2\pi i/h_1}\},\ 
\Spec(T_f) =  \{e^{\nu_q\cdot 2\pi i/h_2}\}
$$ 
then 
$$
\Spec(T_{f\oplus g}) = \{e^{(\mu_p h_2 + \nu_qh_1)\cdot 2\pi i/h_1h_2}\}
$$

\subsection{Simple singularities}
\label{sub24}

 Cf. \cite{AGV} (a), 15.1. They are:
 $$
 x^{n+1},\ n\geq 1,
 \eqno{(A_n)}
 $$
 $$
 x^2y + y^{n-1},\ n\geq 4
 \eqno{(D_n)}
 $$
 $$
 x^4 + y^3
 \eqno{(E_6)}
 $$
 $$
 xy^3 + x^3
 \eqno{(E_7)}
 $$
 $$
 x^5 + y^3
 \eqno{(E_8)}
 $$
 
 Their names come from the following facts: 
 
 --- their lattices of vanishing cycles may be identified with the corresponding root lattices;
 
 --- the monodromy group is identified with the corresponding Weyl group; 
 
 --- the classical monodromy $T_f$ is a Coxeter element, therefore its order
 $h$  is equal to the Coxeter number, and 
 $$
 \Spec(T_f) = \{e^{2\pi i k_1/h},\ldots, e^{2\pi i k_r/h}\}
 $$ 
 where the integers
 $$
 1 = k_1 < k_2 < \ldots < k_r = h - 1,
 $$ 
 are the exponents of our root system.
 
 We will discuss the case of $E_8$ in some details below.   
 
\section{Cartan - Coxeter correspondence}
\label{sec3}

\subsection{Lattices, polarization, Coxeter elements}
\label{sub31}

Let us call {\it a lattice} a pair $(Q, A)$ where $Q$ is a free abelian group, and 
$$
A: Q\times Q\to \mathbb{Z}
$$
a symmetric bilinear map (''Cartan matrix''). We shall identify 
$A$ with a map
$$
A: Q \to Q^\vee := Hom(Q,\mathbb{Z}).
$$
{\it A polarized lattice} is a triple $(Q, A, L)$ where 
$(Q, A)$ is a lattice, and 
$$
L:\ Q\iso Q^\vee
$$
(''variation'', or ''Seifert matrix'') is  an isomorphism such that

\begin{equation}
A = A(L) := L + L^\vee
\end{equation}

where
$$
L^\vee: Q = Q^{\vee\vee}\iso Q^\vee
$$
is the conjugate to $L$.

The {\it Coxeter automorphism} of a polarized lattice is defined by  

\begin{equation}
C = C(L) = - L^{-1}L^\vee \in GL(Q).
\end{equation}

We shall say that the operators $A$ and $C$ are in a {\it Cartan - Coxeter correspondence}.

\vspace{2mm}

\textbf{Example} Let $(Q, A)$ be a lattice, and $\{e_1, \ldots, e_n\}$ 
an ordered $\mathbb{Z}$-base of $Q$. With respect to this base $A$ is expressed as a symmetric matrix $A = (a_{ij}) = A(e_i, e_j)\in \fgl_n(\mathbb{Z})$. Let us suppose that all $a_{ii}$ are even. We define the matrix of $L$ to be the unique upper   
triangular matrix $(\ell_{ij})$ such that $A = L + L^t$ (in particular $\ell_{ii} = a_{ii}/2$; in our examples we will have 
$a_{ii} = 2$.) We will call $L$ the {\it  standard polarization} 
associated to an ordered base. $\square$

\vspace{2mm}

Polarized lattices form a groupoid: 

an isomorphosm  of polarized lattices 
$f:\ (Q_1, A_1, L_1) \iso (Q_2, A_2, L_2)$ is by definition an isomorphism of abelian groups $f: Q_1\iso Q_2$ such that 
$$
L_1(x, y) = L_2(f(x), f(y))
$$ 
(and whence $A_1(x, y) = A_2(f(x), f(y))$). 

\subsection{Orthogonality}
\label{sub32}

\begin{lemma} (i) (orthogonality) {\it 
	$$
	A(x,y) = A(Cx, Cy).
	$$}

(ii) (gauge transformations) {\it For any $P\in GL(Q)$  
	$$
	A(P^\vee L P) = P^\vee A(L)P,\ C(P^\vee L P) = P^{-1}C(L)P.
	$$}
\end{lemma}

$\square$

\subsection{Black/white decomposition and a Steinberg's theorem}
\label{sub34}

Cf. \cite{Stein}, \cite{C}.  Let $\alpha_1, \ldots, \alpha_r$ 
be a base of simple roots of a finite reduced irreducible root 
system $R$ (not necessarily  simply laced). 

Let 
$$
A = (a_{ij}) = (\langle \alpha_i, \alpha_j^\vee \rangle)
$$
be the Cartan matrix. 

Choose a black/white coloring of the set of vertices of the corresponding Dynkin graph $\Gamma(R)$ in such a way that 
any two neighbouring vertices have different colours; this is possible since $\Gamma(R)$ is a tree (cf. \ref{sub-tree}).   

Let us choose an ordering of simple roots in such a way that the first $p$ roots are black, and the last $r - p$ roots are white. 
In this base $A$ has a block form
$$
A = \left(\begin{matrix} 2I_p & X\\ Y & 2I_{r-p}
\end{matrix}\right)
$$
Consider a Coxeter element

\begin{equation}
C = s_1s_2\ldots s_r = C_B C_W,
\end{equation}

where 
$$
C_B = \prod_{i=1}^p s_i,\ C_W = \prod_{i=p+1}^r s_i.
$$
Here $s_i$ denotes the simple reflection corresponding to the root 
$\alpha_i$.

The matrices of $ C_B, C_W$ with respect to the base $\{\alpha_i\}$ 
are
$$  
C_B = \left(\begin{matrix} -I & -X\\ 0 & I\end{matrix}\right),
C_W = \left(\begin{matrix} I & 0\\ -Y & -I\end{matrix}\right), 
$$
so that 

\begin{equation}
C_B + C_W = 2I - A.
\end{equation}

This is an observation due to  R.Steinberg, cf. \cite{Stein}, p. 591. 

\bigskip

We can also rewrite this as follows. 
Set 
$$
L = \left(\begin{matrix} I & 0\\ Y & I\end{matrix}\right),\
U = \left(\begin{matrix} I & X\\ 0 & I\end{matrix}\right). 
$$
Then $A = L + U$, and one checks easily that 

\begin{equation}
C = - U^{-1}L,
\end{equation}

so we are in the situation \ref{sub31}. This explains the name ''Cartan - Coxeter coresspondence''.

\subsection{Eigenvectors' correspondence}
\label{sub35}

\begin{theorem}\label{white-bl} Let 
$$
L = \left(\begin{matrix} I_p & 0\\ Y & I_{r-p}\end{matrix}\right),\   
U = \left(\begin{matrix} I_p & X\\ 0 & I_{r-p}\end{matrix}\right)\ 
$$
be block matrices. Set 
$$
A = L + U,\ C = - U^{-1}L.
$$ 
Let $\mu\neq 0$ be a complex number, $\sqrt{\mu}$ be any of its square roots, and 

\begin{equation}\label{lambda}
\lambda = 2 - \sqrt{\mu} - 1/\sqrt{\mu}.
\end{equation}

Then a vector $v_C = \left(\begin{matrix} v_1\\ v_2\end{matrix}\right)$ 
is an eigenvector of $C$ with  eigenvalue $\mu$ if and only if   
$$
v_A = \left(\begin{matrix} v_1\\ \sqrt{\mu}v_2\end{matrix}\right)
$$
is an eigenvector of $A$ with the eigenvalue $\lambda$\footnote{this formulation has been suggested by A.Givental.}.

\end{theorem}

{\bf Proof}: a direct check. $\square$

\bigskip

\subsubsection{Remark}
Note that the formula (\ref{lambda}) gives two possible values of $\lambda$ corresponding 
to $\pm \sqrt{\mu}$. On the other hand, $\lambda$ does not change 
if we replace $ \mu$ by $\mu^{-1}$. 

In the simplest case of $2\times 2$ matrices the eigenvalues of $A$ are 
$2 \pm (\sqrt{\mu} + \sqrt{\mu^{-1}})$, whereas the eigenvalues of 
$C$ are $\mu^{\pm 1}$.  

\bigskip


\begin{corollary}\label{vp}
	{\it In the notations of} \ref{sub31}, {\it a vector
		$$
		x = \sum x_j\alpha_j
		$$
		is an eigenvector of $A$ with the eigenvalue $2(1 - \cos\theta)$ iff 
		the vector
		$$
		x_c := \sum e^{\pm i\theta/2}x_j \alpha_j
		$$
		where the sign in $e^{\pm i\theta/2}$ is plus if $i$ is a white vertex, and minus otherwise, is an eigenvector of $C$ with eigenvalue $e^{2i\theta}$.}
\end{corollary}

Cf.  \cite{F}.

\begin{proof}
Without loss of generality, we can suppose that $A$ is expressed in a basis of simple roots such that the first $r-p$ ones are white, and the last $p$ roots are black. 

Then $A$ has a block form 

\[A = \begin{pmatrix}
2 I_{r-p} & X \\
Y & 2 I_{p}
\end{pmatrix} = \begin{pmatrix}
I_{r-p} & 0 \\
Y & I_{p}
\end{pmatrix} + \begin{pmatrix}
I_{r-p} & X \\
0 & I_{p} 
\end{pmatrix} = L + U \]

Applying Theorem 1 with 

\[ v_1 = \begin{pmatrix}
e^{i\theta/2} x_1 \\
.. \\
e^{i\theta/2} x_{r-p} \end{pmatrix} \text{ and } v_2 = \begin{pmatrix} e^{-i\theta/2} x_{r-p+1} \\
.. \\
e^{-i\theta/2} x_r \end{pmatrix} \]	

and the well-known eigenvalues of the Cartan matrix $A$, 
\[\lambda = 2 - 2 \cos \theta_k , \text{ with }  \theta_k = 2\pi k/h , k\in \Exp(R) \]
we obtain : $x_c := \sum e^{\pm i\theta/2}x_j \alpha_j $ is an eigenvector of $C$ with the eigenvalue $e^{2 i \theta_k}$ iff $e^{i\theta_k}x = e^{i\theta_k} \sum  x_j \alpha_j$ is an eigenvector of $A$ with the eigenvalue $2 - 2 \cos \theta_k $. $\square$ 
\end{proof}

\subsection{Example: the root systems $A_n$.}
\label{sub41}

We  consider the Dynkin graph of $A_n$ 
with the obvious numbering of the vertices. 

The Coxeter number $h = n + 1$, the set of exponents:
$$
\Exp(A_n) = \{1, 2, \ldots, n\}
$$ 

The eigenvalues of any Coxeter element are $e^{i\theta_k}$, and 
the eigenvalues of the Cartan 
matrix $A(A_n)$ are $2 - 2\cos\theta_k$, $\theta_k = 2\pi k/h$, 
$k\in \Exp(A_n)$.

An eigenvector of $A(A_n)$ with the eigenvalue $2 - 2\cos\theta$ 
has the form 

\begin{equation}
x(\theta) =  (\sum_{k=0}^{n-1} e^{i(n-1 - 2k)\theta},
\sum_{k=0}^{n-2} e^{i(n-2 - 2k)\theta} , \ldots, 
1)
\end{equation}

Denote by $C(A_n)$ the Coxeter element 
$$
C(A_n) = s_1s_2\ldots s_n
$$
Its eigenvector with the eigenvalue $e^{2i\theta}$ is:
$$
X_{C(A_n)} = (\sum_{k=0}^{n-j} e^{2ik\theta})_{1\leq j\leq n}
$$
For example, for $n = 4$:

\[ C_{A_4} = \begin{pmatrix}
0 & 0 & 0 & -1 \\
1 & 0 & 0  &-1 \\
0 & 1 & 0 & -1 \\
0 & 0 & 1 & -1 \end{pmatrix} \text{ and } X_{C(A_4 )} = \begin{pmatrix} 1 + e^{2i\theta} + e^{4i\theta} + e^{6i\theta} \\ 1 + e^{2i\theta} + e^{4i\theta}  \\ 1 + e^{2i\theta} \\ 1 \end{pmatrix}\]

is an eigenvector with eigenvalue  $ e^{2i\theta}$. 

Similarly, for $n = 2$:

\[ C_{A_2} = \begin{pmatrix}
0 & -1 \\
1 & -1 \end{pmatrix},\ X_{C(A_2 )}  = \begin{pmatrix} 1+ e^{2i\gamma} \\ 1 \end{pmatrix}  \]

$\square$

\section{Sebastiani - Thom product; factorization of $E_8$ and $E_6$}
\label{sec4}

\subsection{Join product}
\label{sub33}
Suppose we are given two polarized lattices 
$(Q_i, A_i, L_i)$, $i = 1, 2$.  

Set $Q = Q_1\otimes Q_2$, whence 
$$
L:= L_1\otimes L_2: Q\iso Q^\vee,
$$
and define
$$
A: = A_1*A_2 := L + L^\vee: Q\iso Q^\vee  
$$
The triple $(Q, A, L)$ will be called the {\bf join}, or {\bf Sebastiani - Thom}, product of the polarized lattices 
$Q_1$ and $Q_2$, and denoted by $Q_1*Q_2$. 

Obviously
$$
C(L) = - C(L_1)\otimes C(L_2)\in GL(Q_1\otimes Q_2).
$$

It follows that if $\Spec(C(L_i)) = \{e^{2\pi i k_{i}/h_i},\ k_i\in K_i\}$ then

\begin{equation}
\Spec(C(L)) = \{ - e^{2\pi i(k_{1}/h_1 + k_{2}/h_2)},\ (k_1,k_2)\in K_1\times K_2\}
\end{equation}

\subsection{$E_8$ versus $A_4*A_2*A_1$: elementary analysis}
\label{sub42}

The ranks:
$$
r(E_8) = 8 = r(A_4)r(A_2)r(A_1);
$$
the Coxeter numbers: 
$$
h(E_8) = h(A_4)h(A_2)h(A_1) = 5\cdot 3\cdot 2 = 30.
$$
It follows that
$$
|R(E_8)| = 240 = |R(A_4)||R(A_2)||R(A_1)|.
$$ 

The exponents of $E_8$ are:
$$
1, 7, 13, 19, 11, 17, 23, 29.
$$
All these numbers, except $1$, are primes, and these are all primes $\leq 30$, not dividing $30$. 

They may be determined from the formula 
$$
\frac{i}{5} + \frac{j}{3} + \frac{1}{2} = \frac{30 + k(i,j)}{30},\ 
1\leq i\leq 4,\ 1\leq j\leq 2,
$$
so
$$
k(i, 1)= 1 + 6(i-1) = 1, 7, 13, 19;\ 
$$
$$
k(i,2) = 1 + 10 + 6(i-1) = 11, 17, 23, 29.  
$$

This shows that the exponents of $E_8$ are the same as the exponents 
of 
\newline $A_4*A_2*A_1$.

The following theorem is more delicate. 

\subsection{Decomposition of $Q(E_8)$}
\label{sub43}

\begin{theorem}\label{gab} (Gabrielov, cf. \cite{Gab}, Section 6, Example 3).  There exists a polarization 
	of the root lattice $Q(E_8)$ and an isomorphism 
	of polarized lattices
	
	\begin{equation}\label{gam}
		\Gamma: Q(A_4)*Q(A_2)*Q(A_1) \iso Q(E_8).
	\end{equation}
\end{theorem}

In the left hand side $Q(A_n)$ means the root lattice of $A_n$ with 
the standard Cartan matrix  and the standard polarization
$$
A(A_n) = L(A_n) + L(A_n)^t
$$
where the Seifert matrix $L(A_n)$ is upper triangular.

In the process of the proof, given in \S \ref{sub44} - \ref{sub46} below, the isomorphism $\Gamma$ will be written down explicitly. 

\subsection{Beginning of the proof}
\label{sub44}

For $n = 4, 2, 1$,  
we consider the bases of simple roots $e_1,\ldots, e_n$ in $Q(A_n)$, with scalar products  
given by the Cartan matrices $A(A_n)$. 

The tensor product of three lattices  
$$
Q_* = Q(A_4)\otimes  Q(A_2)\otimes  Q(A_1) 
$$
will be equipped  with the ''factorizable'' basis in the lexicographic order: 
$$
(f_1,\ldots, f_8) := (e_1\otimes e_1\otimes e_1, e_1\otimes e_2\otimes e_1, 
e_2\otimes e_1\otimes e_1, e_2\otimes e_2\otimes e_1,
$$
$$
e_3\otimes e_1\otimes e_1, e_3\otimes e_2\otimes e_1, 
e_4\otimes e_1\otimes e_1, e_4\otimes e_2\otimes e_1).
$$
Introduce a scalar product $(x, y)$ on $Q_*$ given, in the basis $\{f_i\}$, by the matrix
$$
A_* = A_4*A_2*A_1.
$$



\subsection{Gabrielov - Picard - Lefschetz transformations $\alpha_m, 
	\beta_m$}
\label{sub45}

Let $(Q, ( , ))$ be a lattice of rank $r$. We introduce the following 
two sets of transformations $\{\alpha_m\}, \{\beta_m\}$ on the set $Bases-cycl(Q)$ of cyclically ordered bases of $Q$.

If $x = (x_i)_{i\in \mathbb{Z}/r\mathbb{Z}}$ is a base, and $m\in \mathbb{Z}/r\mathbb{Z}$, 
we set

$$
(\alpha_m(x))_i = \left\lbrace \begin{matrix}
x_{m+1} + (x_{m+1}, x_{m})x_m\ & \text{if}\ i = m\\
x_m & \text{if}\ i = m + 1 \\
x_i & \text{otherwise}
\end{matrix}\right.
$$
and 
$$
(\beta_m(x))_i = \left\lbrace \begin{matrix}
x_m & \text{if}\ i = m - 1 \\
x_{m-1} + (x_{m-1}, x_{m})x_m\ & \text{if}\ i = m\\
x_i & \text{otherwise}
\end{matrix}\right.
$$
We define also a transformation $\gamma_m$ by
$$
(\gamma_m(x))_i = \left\lbrace \begin{matrix}
- x_m & \text{if}\ i = m \\
x_i & \text{otherwise}
\end{matrix}\right.
$$  

\subsection{Passage from $A_4*A_2*A_1$ to $E_8$}
\label{sub46}

Consider the base 
$ f = \{f_1, \ldots f_8\}$ of the lattice $Q_* := Q(A_4)\otimes Q(A_2) \otimes Q(A_1)$ described in \S \ref{sub44}, and apply to it the following transformation

\begin{equation}\label{g'}
G' = \gamma_2\gamma_1\beta_4\beta_3\alpha_3\alpha_4\beta_4\alpha_5\alpha_6
\alpha_7\alpha_1\alpha_2\alpha_3\alpha_4\beta_6\beta_3\alpha_1,
\end{equation}
cf. \cite{Gab}, Example 3. Note that
\begin{equation}\label{eqn:gamma}
\gamma_2\gamma_1 = \alpha_1^6,
\end{equation}
cf. \cite{Br}.

Then the base $G'(f)$ has the intersection matrix given by the Dynkin 
graph of $E_8$, with the ordering indicated in Figure \ref{fig1} below.

\begin{figure}
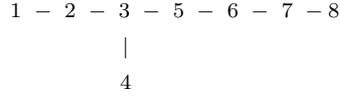

	$$
	1 \  -\   2\    - \   3 \   - \   5 \   - \   6 \   - \   7\ - 8
	$$
	$$
	|\ \ \ \ \ \ \ \ \ \  \ \ \     
	$$
	$$
	4\ \ \ \ \ \ \ \ \ \  \ \ \     
	$$
	
	\caption{Gabrielov's ordering of $E_8$.}
	\label{fig1}       
\end{figure}

This concludes the proof of Theorem \ref{gab} $\square$

\subsection{The induced map of root sets}
\label{sub47}

By definition, the isomorphism of lattices $\Gamma$, (\ref{gam}), induces 
a bijection between the bases
$$
g:\ \{f_1,\ldots, f_8\} \iso \{\alpha_1,\ldots, \alpha_8\}\subset R(E_8).
$$
where in the right hand side we have the base of simple roots, 
and a map
$$
G:\ R(A_4)\times R(A_2)\times R(A_1)\to R(E_8),\ 
G(x,y,z) = \Gamma(x\otimes y \otimes z)  
$$
of sets of the same cardinality $240$ which is not a bijection however: its image consists of $60$ elements. 

Note that the set of vectors $\alpha\in Q(E_8)$ with 
$(\alpha, \alpha) = 2$ coincides with the root system $R(E_8)$, 
cf. \cite{Serre}, Premi\`ere Partie, Ch. 5, 1.4.3.  

\subsection{Passage to Bourbaki ordering}
\label{sub48}

The isomorphism $G'$ (\ref{g'}) is given by a matrix 
$G'\in GL_8(\mathbb{Z})$ such that 
$$
A_G(E_8) = G^{\prime t}A_*G'
$$
where we denoted
$$
A_* = A(A_4)*A(A_2)*A(A_1),
$$
the factorized Cartan matrix, 
and $A_G$ denotes the Cartan matrix of $E_8$ with respect to the 
numbering of roots indicated on Figure \ref{fig1}. 

Now let us pass to the numbering  of vertices of the Dynkin graph 
of type $E_8$ indicated in \cite{B} (the difference with Gabrielov's 
numeration is in three vertices $2, 3$, and $4$).

\begin{figure}[h]
$$
1 \  -\   3\    - \   4 \   - \   5 \   - \   6 \   - \   7\ - 8
$$
$$
|\ \ \ \ \ \ \ \ \ \  \ \ \     
$$
$$
2\ \ \ \ \ \ \ \ \ \  \ \ \     
$$	
	
	\caption{Bourbaki ordering of $E_8$.}
	\label{fig2}    
\end{figure}

The Gabrielov's Coxeter element (the full monodromy) in the Bourbaki numbering looks as follows: 
$$
C_{G}(E_8 ) = s_1 \circ s_3 \circ s_4 \circ s_2 \circ s_5 \circ s_6 \circ s_7 \circ s_8
$$


\begin{lemma}{\it Let $A(E_8)$ be the standard Cartan matrix of 
	$E_8$ from} [B]:
\[ A(E_8) = \begin{pmatrix}
2 & 0 & -1 & 0 & 0 & 0 & 0 & 0 \\
0 & 2 & 0 & -1 & 0 & 0 & 0 & 0 \\
-1 & 0 & 2 & -1 & 0 & 0 & 0 & 0 \\
0 & -1 & -1 & 2 & -1 & 0 & 0 & 0 \\
0 & 0 & 0 & -1 & 2 & -1 & 0 & 0 \\
0 & 0 & 0 & 0 & -1 & 2 & -1 & 0 \\
0 & 0 & 0 & 0 & 0 & -1 & 2 & -1 \\
0 & 0 & 0 & 0 & 0 & 0 & -1 & 2 \end{pmatrix}. \]
{\it Then
	$$
	A(E_8) = G^{ t}A_*G
	$$
	and
	$$
	C_G(E_8) = G^{-1}C_* G
	$$
	where
	$$
	C_* = C(Q(A_4)*Q(A_2)*Q(A_1)) = C(A_4)\otimes C(A_2)\otimes C(A_1),
	$$
	is  the factorized Coxeter element, 
	and   
	\[ G = \begin{pmatrix}
	0 & 0 & 0 & 1 & -1 & 0 & 0 & 0 \\
	-1 & 1 & 0 & 0 & 0 & 0 & 0 & 0 \\
	0 & 0 & -1 & 1 & 0 & 0 & 0 & 0 \\
	-1 & 1 & -1 & 0 & 0 & 1 & 0 & 0 \\
	0 & 1 & -1 & 0 & 0 & 0 & 1 & 0 \\
	-1 & 1 & -1 & 0 & 0 & 0 & 1 & 0 \\
	0 & 1 & -1 & 0 & 0 & 0 & 0 & 1 \\
	0 & 1 & -1 & 0 & 0 & 0 & 0 & 0 \end{pmatrix} 
	\eqno{(3.8.1)} \]
	Here
	$$
	G = G'P
	$$
	where $P$ is the permutation matrix of passage from the Gabrielov's ordering in Figure \ref{fig1} to the Bourbaki ordering in Figure \ref{fig2}}
\end{lemma}

\subsection{Cartan eigenvectors of $E_8$}
\label{sub491}

To obtain the Cartan eigenvectors of $E_8$, one should pass from $C_G(E_8)$ to 
the ''black/white'' Coxeter element (as in \S \ref{sub34})   

\[C_{BW}(E_8) = s_1 \circ s_4 \circ s_6 \circ s_8 \circ s_2 \circ s_3 \circ s_5 \circ s_7 \]

Any two Coxeter elements are conjugate in the Weyl group $W(E_8)$.

The elements $C_G(E_8)$ and $C_{BW}(E_8)$  
are conjugate by the following element of  $W(E_8)$:

\[ C_G(E_8) = w^{-1}C_{BW}(E_8)w \]

where

\[ w = s_7 \circ s_5 \circ s_3 \circ s_2 \circ s_6 \circ s_4 \circ s_5 \circ s_1 \circ s_3 \circ s_2 \circ s_4 \circ s_1 \circ s_3 \circ s_2 \circ s_1 \circ s_2  \]

This expression for $w$ can be obtained using an algorithm described in \cite{C}, 
cf. also \cite{Br}.  

Thus, if $x_*$ is an eigenvector of $C_*(E_8)$ then 

\[ x_{BW} = wG^{-1}x_* \]

is an eigenvector of $C_{BW}(E_8)$. But we know the eigenvectors of 
$C_*(E_8)$, they are all factorizable. 

This provides the eigenvectors of $C_{BW}(E_8)$, which in turn 
have very simple relation to the eigenvectors of $A(E_8)$, 
due to Theorem \ref{vp}.

\bigskip

{\bf Conclusion: an expression for the eigenvectors of $A(E_8)$.}

\bigskip 

Let $\theta = \frac{a\pi}{5},\ 1\leq a\leq 4$, $\gamma = \frac{b\pi}{3},\ 1\leq b\leq 2$, $\delta = \frac{\pi}{2}$,
$$
\alpha = \theta + \gamma + \delta = \pi + \frac{k\pi}{30},
$$
$$
k\in \{1, 7, 11, 13, 17, 19, 23, 29\}.
$$
The $8$ eigenvalues of $A(E_8)$ have the form
$$
\lambda(\alpha) = \lambda(\theta,\gamma) = 2 - 2\cos \alpha
$$ 
An eigenvector of $A(E_8)$ with the eigenvalue 
$\lambda(\theta,\gamma)$ is 

\[ X_{E_8}(\theta,\gamma) = \begin{pmatrix} 
\cos(\gamma + \theta - \delta) + \cos ( \gamma - 3 \theta - \delta) + \cos (\gamma - \theta - \delta) \\
\cos( 2 \gamma + 2 \theta)\\
\cos(2\gamma) + \cos (2\gamma + 2 \theta) + \cos ( 2\gamma - 2 \theta) + \cos (4\theta) + \cos(2 \theta) \\
\cos(\gamma + 3 \theta - \delta) + \cos ( \gamma + \theta - \delta ) + \cos ( -\gamma + 3 \theta - \delta) \\
2 \cos(2\gamma) + 2 \cos(2 \gamma + 2 \theta) + \cos(2\gamma - 2 \theta) + \cos(2 \gamma + 4 \theta) + \cos(4\theta) + 2 \cos (2 \theta) + 1 \\
\cos(\gamma + 3 \theta - \delta) + \cos( \gamma + \theta - \delta) \\
\cos (2 \gamma) + \cos ( 2 \theta - 2 \delta) \\
\cos (\gamma - \theta - \delta ) \end{pmatrix} \]

One can simplify it as follows:

\begin{equation}
 X_{E_8}(\theta,\gamma) =  - \begin{pmatrix} 
 2 \cos (4\theta) \cos(\gamma - \theta - \delta) \\
 - \cos( 2 \gamma + 2 \theta ) \\
 2 \cos^{2} (\theta) \\
 -2 \cos(\gamma) \cos (3\theta - \delta) - \cos ( \gamma + \theta - \delta ) \\
 -2 \cos(2\gamma + 3 \theta) \cos (\theta) + \cos (2 \gamma ) \\
 -2 \cos \theta \cos (\gamma + 2\theta - \delta ) \\
 -2\cos(\gamma + \theta - \delta) \cos(\gamma -\theta + \delta) \\
 - \cos (\gamma - \theta - \delta) \end{pmatrix} 
\end{equation}

\subsection{Perron - Frobenius and all that}
\label{sub492}

The Perron - Frobenius eigenvector corresponds to the eigenvalue

\[ 2 - 2 \cos \frac{\pi}{30},  \]

and may be chosen as

\[ v_{PF} = \begin{pmatrix} 
2 \cos \frac{\pi}{5} \cos  \frac{11\pi}{30}   \\
\cos \frac{\pi}{15}  \\
2\cos^{2}\frac{\pi}{5} \\ 
2\cos  \frac{2\pi}{30} \cos \frac{\pi}{30} \\      
2 \cos \frac{4\pi}{15}  \cos \frac{\pi}{5} + \frac{1}{2} \\  
2 \cos \frac{\pi}{5} \cos \frac{7\pi}{30} \\
2 \cos \frac{\pi}{30} \cos \frac{11\pi}{30}\\  
\cos \frac{11\pi}{30} \end{pmatrix} \] 

Ordering its coordinates in the increasing order, we obtain

\[ v_{PF <} = \begin{pmatrix} 
\cos \frac{11\pi}{30} \\
2 \cos \frac{\pi}{5} \cos  \frac{11\pi}{30}  \\
2 \cos \frac{\pi}{30} \cos \frac{11\pi}{30} \\
\cos \frac{\pi}{15} \\
2 \cos \frac{\pi}{5} \cos \frac{7\pi}{30} \\
2\cos^{2}\frac{\pi}{5}  \\
2 \cos \frac{4\pi}{15}  \cos \frac{\pi}{5} + \frac{1}{2} \\
2\cos  \frac{2\pi}{30} \cos \frac{\pi}{30} \end{pmatrix} \]

In the Ref. \cite{Z}, A. B. Zamolodchikov obtains 
the following expression for the PF vector: 

\[ v_{Zam}(m) = \begin{pmatrix}
m \\
2 m \cos \frac{\pi}{5} \\
2 m \cos \frac{\pi}{30} \\
4 m \cos \frac{\pi}{5} \cos \frac{7\pi}{30} \\
4 m \cos \frac{\pi}{5} \cos \frac{2\pi}{15} \\
4 m \cos \frac{\pi}{5} \cos \frac{\pi}{30} \\
8 m \cos^{2} \frac{\pi}{5} \cos \frac{7\pi}{30} \\
8 m \cos^{2} \frac{\pi}{5} \cos \frac{2\pi}{15} \end{pmatrix} \] 

Setting $m = \cos \frac{11\pi}{30}$, we find indeed : 

\[ v_{PF <} =  v_{Zam} (\cos \frac{11\pi}{30} )  \]




\subsection{Factorization of $E_6$}
\label{sub-E6}

\begin{theorem}
	(Gabrielov, cf. \cite{Gab}, Section 6, Example 2).  There exists a polarization 
	of the root lattice $Q(E_6)$ and an isomorphism 
	of polarized lattices
	
	\begin{equation}\label{gam}
	\Gamma_{E_6} : Q(A_3)*Q(A_2)*Q(A_1) \iso Q(E_6).
	\end{equation}
\end{theorem}

The  proof is exactly the same as for $Q(E_8)$. The passage from $A_3 * A_2 * A_1$ to $E_6$ is obtained by the following transformation

\[ G'_{E_6} = \gamma_4 \gamma_1 \alpha_1 \alpha_2 \alpha_3 \alpha_4 \beta_6 \beta_3 \alpha_1 \]

cf. \cite{Gab}, Example 2.

After a passage from Gabrielov's ordering to Bourbaki's, we obtain a transformation 

\[ G_{E_6} = \begin{pmatrix}
0&-1&1&0&0&0\\
-1&0&1&0&0&0\\
0&-1&0&1&0&0\\
-1&0&0&0&1&0\\
0&0&0&0&0&1\\
-1&0&0&0&0&1
\end{pmatrix} \in GL_{6} ( \mathbb{Z} )\]

such that 
\[ A(E_6) = G_{E_6}^{t} A_{*} G_{E_6} \text{ and }C_{G} (E_6) = G_{E_6}^{-1} C_{*} G_{E_6} \]  
where $ A_{*} = A(A_3) * A(A_2) * A(A_1)$ and $ C_{*} = C(A_3) \otimes C(A_2) \otimes C(A_1)$ and 
\[ C_{G}(E_6) = s_1 \circ s_3 \circ s_4 \circ s_2 \circ s_5 \circ s_6  \]
 $C_{G}(E_6)$ is the Gabrielov's Coxeter element in the Bourbaki numbering, cf. \cite{B}.  
 
 Let $C_{BW} (E_6) = s_1 \circ s_4 \circ s_6 \circ s_2 \circ s_3 \circ s_5$ be the "black/white" Coxeter element. $C_{G}(E_6)$ and $C_{BW}(E_6)$ are conjugated by the following element of the Weyl group $W(E_6)$ :
 
 \[ v = s_5 \circ s_3 \circ s_2 \circ s_4 \circ s_1 \circ s_3 \circ s_3 \circ s_1 \circ s_2 \]

Thus, if $x_{*}$ is an eigenvector of $C_{*} ( E_6)$ then $x_{BW} = v G^{-1}_{E_6} x_{*}$ is an eigenvector of $C_{BW}(E_6)$. 

Finally, let $\theta = \frac{a\pi}{4} , 1 \leq a \leq 3$, $\gamma = \frac{b\pi}{3}, 1 \leq b \leq 2$, $\delta = \frac{\pi}{2}$ and

\[ \alpha = \theta + \gamma + \delta \]

The 6 eigenvalues of $A(E_6)$ have the form $\lambda (\alpha) = \lambda(\theta, \gamma) = 2 - 2\cos \alpha$. An eigenvector of $A(E_6)$ with the eigenvalue $\lambda (\alpha)$ is

\[ X_{E_6} ( \theta, \lambda) = \begin{pmatrix}
\cos \left(3\gamma + 3 \theta - \delta \right) \\
2 \cos^{2} \theta \\
-2 \cos \left( 3 \gamma + 3 \theta - \delta\right) \cos \left(\gamma + \theta - \delta\right) \\
-4 \cos^{2} \theta \cos \left( \gamma + \theta - \delta\right) \\
1 - 2 \cos \left( 2 \gamma + 3 \theta\right) \cos \theta \\
-2 \cos(\gamma) \cos \left(\theta - \delta\right)
\end{pmatrix}    \]

\section{Givental's $q$-deformations}
\label{sec5}

\subsection{$q$-deformations of Cartan matrices}
\label{sub-cartan}

Let $A = (a_{ij})$ be a $n\times n$ complex matrix. We will say that $A$ is 
a {\it generalized Cartan matrix} if

\bigskip

(i) for all $i\neq j$,  $a_{ij}\neq 0$ implies $a_{ji}\neq 0$;

(ii) all $a_{ii} = 2$.

\bigskip

If only (i) is fulfilled, we will say that $A$ is 
a {\it pseudo-Cartan matrix}.  

\bigskip

We associate to a pseudo-Cartan matrix $A$ an unoriented graph $\Gamma(A)$ with vertices 
$1, \ldots, n$, two vertices $i$ and $j$ being connected by an edge  $e = (ij)$ iff $a_{ij}\neq 0$. 

Let $A$ be a generalized Cartan matrix. 
There is a unique decomposition 
$$
A = L + U
$$
where $L = (\ell_{ij})$ (resp. $U = (u_{ij})$) is lower (resp. upper) triangular, with $1$'s 
on the diagonal. 

We define a $q$-deformed Cartan matrix by 
	$$
	A(q) = qL + U
	$$
	
This definition is inspired by the $q$-deformed Picard - Lefschetz 
theory developed by Givental, \cite{Giv}.

\begin{theorem}\label{q-thm} Let $A$ be a  generalized Cartan matrix 
such that $\Gamma(A)$ is a tree.  

(i) The eigenvalues of $A(q)$ have the form
\begin{equation}\label{q-lambda}
\lambda(q) = 1 + (\lambda - 2)\sqrt{q} + q
\end{equation}
where $\lambda$ is an eigenvalue of $A$. 

(ii) There exist integers $k_1, \ldots, k_n$  such that if $x = (x_1, \ldots, x_n)$ 
is an eigenvector of $A$ for the eigenvalue $\lambda$ then
\begin{equation}\label{q-vect}
x(q) = (q^{k_1/2}x_1,\ldots, q^{k_n/2}x_n) 
\end{equation}
is an eigenvector of $A(q)$ for the eigenvalue $\lambda(q)$. 
\end{theorem}

The theorem will be proved after some preparations. 

\subsection{}\label{sub-tree} Let $\Gamma$ be an unoriented tree 
with a finite set of vertices $I = V(\Gamma)$.

Let us pick  a root of  $\Gamma$, and partially order its vertices 
by taking the minimal vertex $i_0$ to be the bottom of the root, and then 
going ''upstairs''.  This defines an orientation on $\Gamma$.

\begin{lemma}\label{tree-lemma} Suppose we are given a nonzero complex number $b_{ij}$ for each 
edge $e = (ij), i < j$ of $\Gamma$. 
There exists a collection of nonzero complex numbers $\{c_i\}_{i\in I}$ such that 
$$
b_{ij} = c_j/c_i,\ i < j.
$$
for all edges $(ij)$. 

We can choose the numbers $c_i$ in such a way that they are products 
of some numbers $b_{pq}$.

\end{lemma}


{\bf Proof.} Set $c_{i_0} = 1$ for the unique minimal vertex $i_0$, 
and then define the other $c_i$ one by one, by going upstairs, 
and using as a definition
$$
c_j := b_{ij}c_i,\ i < j.
$$
Obviously, the numbers $c_i$ defined in such a way, are products 
of $b_{pq}$. $\square$ 

\begin{lemma}\label{two-matrices} Let $A = (a_{ij})$ and $A' = (a'_{ij})$ be two 
pseudo-Cartan matrices with $\Gamma(A) = \Gamma(A')$. Set 
$b_{ij}:= a'_{ij}/a_{ij}$. Suppose  that
\begin{equation}\label{eq-two}
b_{ij} = b_{ji}^{-1}. 
\end{equation} 

for all $i\neq j$, and $a_{ii} = a'_{ii}$ for all $i$. Then there exists a diagonal matrix 
$$
D = \Diag(c_1,\ldots, c_r)
$$
such that $A' = D^{-1}AD$.

Moreover, the numbers $c_i$ may be chosen to be products of some 
$b_{pq}$.
\end{lemma}

{\bf Proof.} Let us choose a partial order $<_p$ on the set of vertices $V(\Gamma)$ as in \ref{sub-tree}. 

{\it Warning.} This partial order {\it differs} in general from 
the standard total order on $\{1,\ldots, n\}$. 

Let us apply Lemma\ \ref{tree-lemma} to the collection of numbers 
$\{b_{ij},\ i<_p j\}$. We get a sequence of numbers $c_{ij}$ such 
that 
$$
b_{ij} = c_j/c_i
$$
for all $i <_p j$. The condition (\ref{eq-two}) implies that this 
holds true for all $i\neq j$.

By definition, this is equivalent to 
$$
a'_{ij} = c_i^{-1}a_{ij}c_j,
$$
i.e. to $A' = D^{-1}AD$. $\square$

\subsection{\bf{Proof of Theorem \ref{q-thm}.}}
\label{proof-q-thm} 

Let us consider two matrices: $A(q) = (a(q)_{ij})$ with $a(q)_{ii} = 
1 + q$
$$
a(q)_{ij} = \left\{\begin{matrix} a_{ij} & \text{if\ }i < j\\
                                                     qa_{ij} & \text{if\ }i > j
\end{matrix}\right.
$$
and
$$
A'(q) = \sqrt{q}A + (1 -  \sqrt{q})^2 I = (a(q)'_{ij})
$$
with $a(q)'_{ii} = 1 + q$ and $a(q)'_{ij} = \sqrt{q}a(q)_{ij}$, 
$i\neq j$. 

Thus, we can apply Lemma\ \ref{two-matrices} to  $A(q)$ and   
$A'(q)$. 
So, there exists a diagonal matrix $D$ as above such that 
$$
A(q) = D^{-1}A'(q)D.
$$
But the eigenvalues of $A'(q)$ are obviously 
$$
\lambda(q) =  \sqrt{q}\lambda + (1 -  \sqrt{q})^2  = 
1 + (\lambda - 2)\sqrt{q} + q.
$$
If $v$ is an eigenvector of $A$ for $\lambda$ then $v$ is an eigenvector of $A'(q)$ for $\lambda(q)$, and $Dv$ will be  an eigenvector of $A(q)$ for $\lambda(q)$. $\square$

\subsection{{\bf Remark} (M.Finkelberg)}
\label{sub53}

The expression (\ref{q-lambda}) resembles 
the number of points of an elliptic curve $X$ over a finite field 
$\mathbb{F}_q$. To appreciate better this resemblance, note that in all 
our examples $\lambda$  has the form 
$$
\lambda = 2-2 \cos \theta,
$$
so if we set 
$$
\alpha =  \sqrt{q}e^{i\theta}
$$
(''a Frobenius root'') then $|\alpha| = \sqrt{q}$, and 
$$
\lambda(q) = 1 - \alpha - \bar\alpha + q,
$$ 
cf. \cite{IR}, Chapter 11, \S 1, \cite{Kn}, Chapter 10, Theorem 10.5. 

So, the Coxeter eigenvalues $e^{2i\theta}$ may be seen as analogs of ''Frobenius roots of an elliptic curve over $\mathbb{F}_1$''. 


\subsection{\bf Examples. }
\label{sub54}

\subsubsection{Standard deformation for $A_n$}

Let us consider the 
following $q$-deformation of $A = A(A_n)$: 
$$
A(q) = \left(\begin{matrix} 
1 + q & - 1 & 0 & \ldots & 0\\
- q & 1 + q & - 1 & \ldots & 0\\
\ldots & \ldots & \ldots & \ldots & \ldots \\
0 & \ldots & 0 & - q & 1 + q
\end{matrix}\right)
$$
Then 
$$
\Spec(A(q)) = \{\lambda(q) := 1 + (\lambda - 2)\sqrt{q} + q|\ \lambda\in \Spec(A(1))\}.
$$
If $x = (x_1,\ldots, x_n)$ is an eigenvector of $A = A(1)$ 
with eigenvalue $\lambda$ then
$$
x(q) = (x_1, q^{1/2}x_2,\ldots, q^{(n-1)/2}x_n)
$$
is an eigenvector of $A(q)$ with eigenvalue $\lambda(q)$. 

\subsubsection{Standard deformation for $E_8$}
\label{sub55}

A $q$-deformation: 
\[ A_{E_8}(q) = \begin{pmatrix}
1+q & 0 & -1 & 0 & 0 & 0 & 0 & 0 \\
0 & 1+q & 0 & -1 & 0 & 0 & 0 & 0 \\
-q &0& 1+q & -1 & 0 & 0 & 0 & 0 \\
0&-q&-q& 1+q & -1 & 0 & 0 & 0 \\
0 & 0 & 0 & -q & 1+q & -1 & 0 & 0\\
0& 0 & 0 & 0 & -q & 1+q & -1 & 0 \\
0 & 0 & 0 & 0 & 0 & -q & 1+q & -1 \\
0 & 0 & 0 & 0 & 0 & 0 & -q & 1+q \end{pmatrix} \]

Its eigenvalues are 
$$
\lambda(q) = 
1+q+(\lambda-2)\sqrt{q} = 1+q-2 \sqrt{q} \cos \theta
$$ 
where $\lambda = 2-2 \cos \theta$ is an eigenvalue of  $A(E_8)$. 

If $X = (x_1 , x_2 , x_3 , x_4 , x_5 , x_6 , x_7 , x_8)$ is an eigenvector of  $A(E_8)$ for the eigenvalue $\lambda$, then 

\begin{equation}\label{eq1}
X = ( x_1 , \sqrt{q} x_2 , \sqrt{q} x_3 , q x_4 , q \sqrt{q} x_5 , q^2 x_6 , q^2 \sqrt{q} x_7 , q^3 x_8 )
\end{equation}
is an eigenvector of $A_{E_8}(q)$  for the eigenvalue $\lambda(q)$.

\bigskip



\section{A physicist's appendix: 
cobalt niobate producing an $E_8$ chord}

In this Section, we briefly describe  the relation of Perron-Frobenius components, 
in the case of $R = E_8$, to  the physics of certain magnetic systems as anticipated 
in a pioneering theoretical work \cite{Z} and possibly observed 
in a beautiful neutron scattering experiment \cite{Coldea}. 

\subsection{One-dimensional Ising model in a magnetic field}

(a) {\it The Ising Hamiltonian}

Let $W = \BC^2$. Recall three Hermitian Pauli matrices:  
$$
\sigma^x = \left(\begin{matrix}0 & 1\\1 & 0\end{matrix}\right), \,\,
\sigma^y = \left(\begin{matrix}0 & -i\\i & 0\end{matrix}\right), \,\, 
\sigma^z = \left(\begin{matrix}1 & 0\\0 & -1\end{matrix}\right). 
$$
The $\BC$-span of $\sigma^x,  \sigma^y, \sigma^z$ inside $End(W)$ 
is a complex Lie algebra $\fg = \fsl(2,\BC)$; the $\BR$-span of the 
anti-Hermitian matrices $i\sigma^x,  i\sigma^y, i\sigma^z$ is a real Lie subalgebra $\fk = \fsu(2)\subset \fg$. The resulting representation 
of $\fg$ (or $\fk$) on $W$ is what  
physicists refer to as the ``spin-$\frac{1}{2}$ 
representation''. 

For a natural $N$, consider a $2^N$-dimensional tensor product
$$
V = \otimes_{n=1}^N W_n
$$
with all $W_n = W$. We are interested in the spectrum of the following 
linear operator $H$ acting on $V$: 
\begin{equation}\label{eqn:ising-ham}
H = H(J,h_z, h_x) = - J\sum_{n=1}^{N} \sigma^z_n\sigma^z_{n+1} - 
h_z\sum_{n=1}^{N} \sigma^z_n
- h_x\sum_{n=1}^{N} \sigma^x_n ,
\end{equation} 
where $J, h_x, h_z$ are positive real numbers.  Here for $A:\ W \lra W$,  $A_n: V\lra V$ denotes an operator acting as $A$ 
on the $n$-th tensor factor and as the identity on all the other factors. 
By definition, $A_{N+1} := A_1$. 

In keeping with the conditions of the experiment \cite{Coldea}, everywhere 
below we assume that $N$ is very large ($N >> 1$), and that $0 < h_z << J$.

The space $V$ arises as the space of states of a quantum-mechanical 
model   describing a chain of $N$ atoms on the plane $\BR^2$ with 
coordinates $(x,z)$. The chain is parallel to the $z$ axis, and is subject 
to a magnetic field with a component $h_z$ along the chain, and 
a component $h_x$ along the $x$-axis. The $W_n$ is the space of 
states of the $n$-th atom. Only the nearest-neighbor atoms interact, 
and the $J$ parameterizes the strength of this interaction.

The operator $H$ in the Eq. (\ref{eqn:ising-ham}) is 
called the Hamiltonian, and its eigenvalues $\epsilon$ correspond to the 
energy of the system. It is an Hermitian operator (with respect to an 
obvious Hermitian scalar product on $V$), thus all its eigenvalues are real. 

Consider also the {\it translation} operator $T$, acting as follows: 
\begin{equation}
T(v_1\otimes v_2\otimes\ldots \otimes v_N) = 
v_2\otimes v_3\otimes\ldots \otimes v_N\otimes v_1,
\end{equation} 
The operator $T$ is unitary, and commutes with $H$.

An eigenvector $v_0\in V$ of $H$ with the lowest  
energy eigenvalue $\epsilon_0$ is called the ground state. 

What happens as $h_x$ varies, at fixed $J$ and $h_z$?
When $h_x << J$, the ground state $v_0$ is close to the ground 
state $v_J$ of the operator $H_J = H(J,0,0)$: 
$$
v_J = \otimes_{n=1}^N v_n^z ,
$$
where $v^z_n$ is an eigenvector of $\sigma^z$ in $W_i$ 
with eigenvalue $1$. Thus, the state $v_J$ is 
interpreted as ``all the spins pointing along the $z$-axis''.

In the opposite limit, when $h_x >> J$, 
the ground state $v_0$ is close to the ground state $v_x$ 
of the operator $H_x = H(0,0,h_x)$: 
$$
v_x = \otimes_{n=1}^N v_n^x ,
$$
where $v^x_n$ is an eigenvector of $\sigma^x$ in 
$W_n$ with eigenvalue $1$. Thus, the state $v_x$ is 
interpreted as ``all the spins pointing along the $x$-axis''.

As a function of $h_x$ at fixed $J$ and $h_z$, 
the system has two phases. There is a critical value $h_x = h_c$, 
of the order of $J/2$ : for $h_x < h_c$, the ground state $v_0$ is 
close to $v_J$, and  one says that the chain is in the {\it ferromagnetic} 
phase. By contrast, for $h_x > h_c$, the ground state $v_0$  is close 
to $v_x$, and one says that the chain is in the {\it paramagnetic} phase. 
(The transition between the two phases is far less trivial 
than the spins simply turning to follow the field upon increasing $h_x$: to find 
out more, curious reader is encouraged to consult the Ref. \cite{Chakrabarti}.)

(b) {\it Elementary excitations at $h_x = h_c$}

Zamolodchikov's theory, \cite{Z}, says something 
spectacularly precise about the next few, 
after $\epsilon_0$, eigenvalues (``energy levels'') of a nearly-critical 
Hamiltonian $H_c := H(J,h_z << J, h_c)$. 
To see this, notice that the possible eigenvalues of the translation operator $T$ 
have the form $e^{2\pi i k/N}$, with $-N/2\leq k \leq N/2$; let us call the number 
$$
p =2\pi k /N
$$
the 
{\it momentum} of an eigenstate. 
Since $H$ commutes with $T$,  each eigenspace 
\newline 
$V_\epsilon := \{v\in V|\ H_c v = \epsilon v\}$ decomposes further as per
$$
V_\epsilon  = \oplus_p\ V_{p,\epsilon}, 
$$ 
where  
$$
V_{p,\epsilon} := \{v\in V|\ H_cv =  \epsilon v, 
Tv = e^{i p}v\}
$$ 
Let us add a constant to $H_c$ in such a way that the ground state energy
$\epsilon_0$ becomes $0$ and, on the plane $P$ with coordinates $(p, \epsilon)$, 
let us mark all the points, for which $V_{p,\epsilon} \neq 0$.

Zamolodchikov predicted \cite{Z},  that there exist $8$ numbers 
$0 < m_1 < \ldots < m_8$ with the following property. Let us draw on $P$ eight 
hyperbolae
\begin{equation} 
\label{eq:dispersion}
\Hyp_i: \ \epsilon = \sqrt{m_i^2 + p^2},\ 1\leq i \leq 8.
\end{equation}

All the marked points will be located: 

---- either in a vicinity of one of the hyperbolae 
$\Hyp_i$ (in the limit $N\lra \infty$ they will all lie on these hyperbolae). 

\begin{figure}
\includegraphics[scale=0.3]{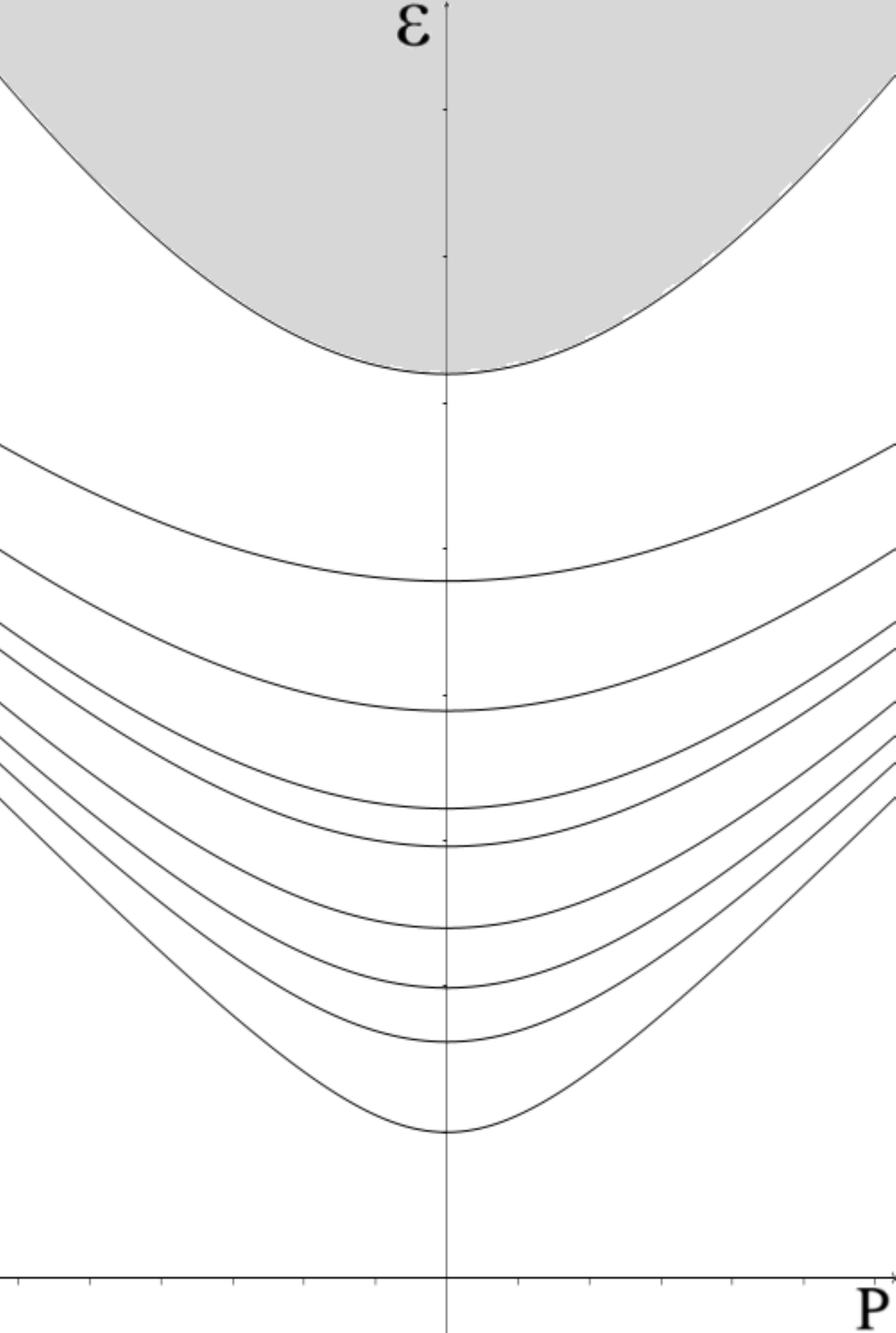}
\caption{The expected joint spectrum of the operators $T, H$.} 
\label{fig:hyperbol}
\end{figure}

---- or in a shaded region separated from these hyperbolae 
as shown in the Fig. \ref{fig:hyperbol}. 

The states $v\in V_{p,\epsilon}$ 
with $(p, \epsilon)\in \Hyp_i$ are called {\it elementary 
excitations}.  The numbers $m_i$ are called their {\it masses}. 


The vector
\begin{equation}  
\label{eq:masses}
\vec{m} = (m_1, \ldots, m_8)
\end{equation} 
is proportional to the Perron - Frobenius $v_{PF<}$ for $E_8$ from \ref{sub492}, 
whose normalized approximate value is
\begin{equation}\label{eq:e8-masses}
v_{PF<} = (1, 1.62, 1.99, 2.40, 2.96, 3.22, 3.89, 4.78)
\end{equation} 

These low-lying excitations (hyperbolae) are 
observable: one may be able to see them  

(a) in a computer simulation, or 

(b) in a neutron scattering experiment.


\subsection{Neutron scattering experiment}

The paper \cite{Coldea} reports the results of a magnetic neutron scattering 
experiment on cobalt niobate CoNb$_2$O$_6$, a material that can be pictured 
as a collection of parallel non-interacting
one-dimensional chains of atoms. We depict such a chain as a straight line, 
parallel to the $z$-axis in our physical space $\BR^3$ with coordinates $x, y, z$. 

The sample, at low temperature $T < 2.95$K (Kelvin), was subject to an 
external magnetic field with components $(h_z,h_z)$, with the $h_x$ at 
the critical value $h_x = h_c$, 
and with $h_z <<  h_c$. 
The system may be described as the Ising chain 
with a nearly-critical Hamiltonian $H = H(J, h_z << h_c, h_c)$ of the Eq. 
(\ref{eqn:ising-ham}). The experiment \cite{Coldea} may be interpreted 
with the help of the following (oversimplified) theoretical picture.

Consider a neutron scattering off the sample. If the incident neutron has  
energy $\epsilon$ and momentum $p$, and scatters off with  energy $\epsilon'$ 
and momentum $p'$, the energy and momentum  conservation laws imply that 
the differences, called  energy and momentum transfers 
$\omega = \epsilon - \epsilon', q = p - p'$, are absorbed by the sample.

The energy transfer cannot be arbirtary. Suppose that, prior to scattering the neutron, 
the sample was in the ground state $v_0$; upon scattering the neutron, it undergoes 
a transition to a state that is a linear combination of the eight elementary 
excitations $v\in V_{p,\epsilon}$.

\begin{figure}
\includegraphics[scale=0.3]{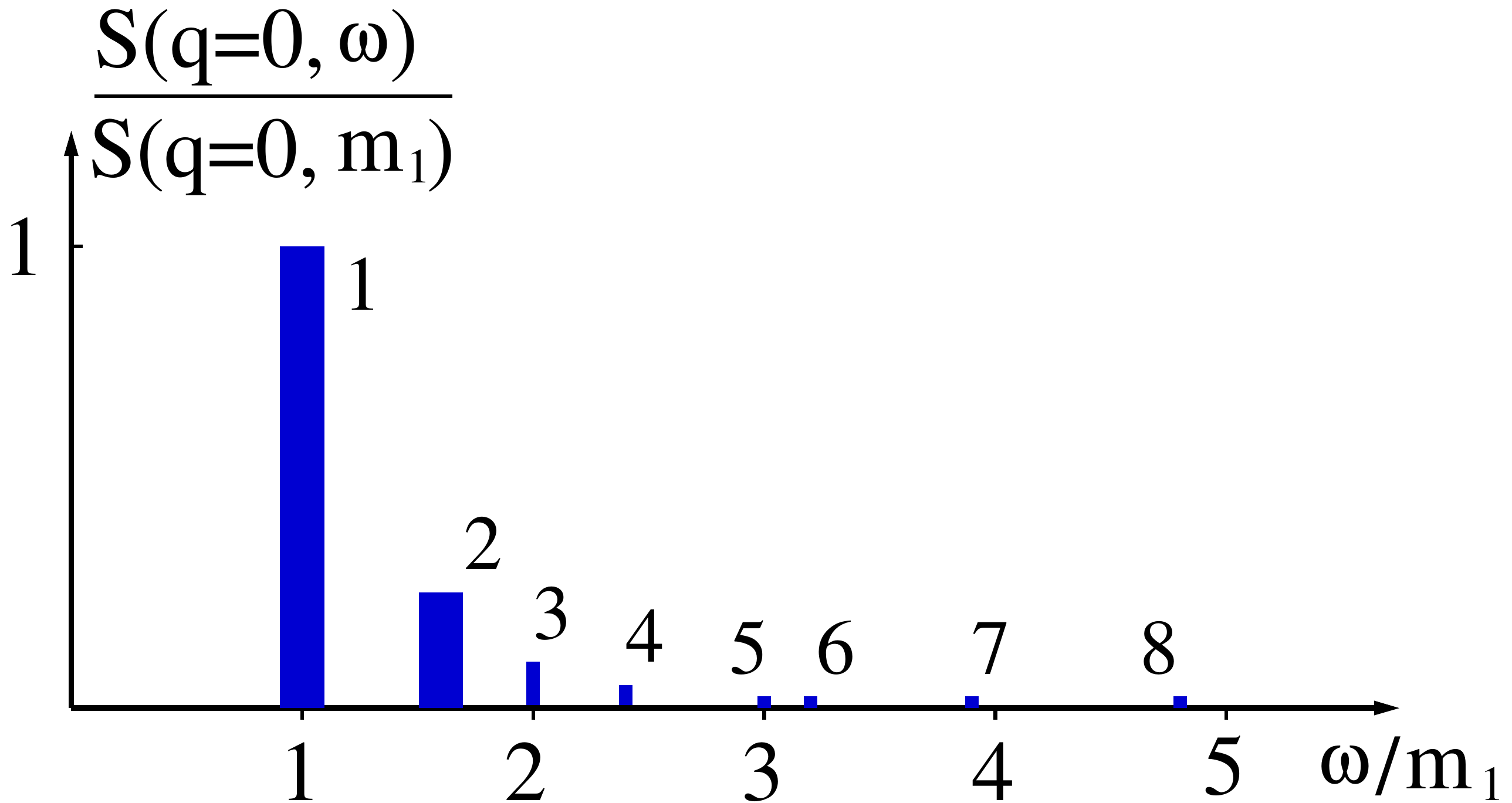}
\caption{A sketch of the  scattering 
intensity $\mathcal{S}(0,\omega)$ at zero momentum peaks 
relative to $\mathcal{S}(0, m_1)$, against the $\omega / m_1$ 
ratio. The two leftmost peaks shown by thick lines correspond 
to the excitations with the masses $m_1$ and $m_2$, 
that were resolved in the experiment \cite{Coldea}. 
The experimentally found mass ratio $m_2/m_1$ is consistent 
with $\frac{m_2}{m_1} = \frac{1 + \sqrt{5}}{2}$, as per the 
expression for the $v_{Zam}(m)$ in the Subsection \ref{sub492}. 
}
\label{fig:mass.peaks}
\end{figure}

We will be interested in neutrons that scatter off with zero 
momentum transfer. The Zamolodchikov theory \cite{Z} predicted, that the neutron scattering intensity $\mathcal{S}(0,\omega)$ should have 
peaks at $\omega = m_a$, ($a = 1, ..., 8$) of the Eq. (\ref{eq:e8-masses}). 
At zero momentum transfer, a neutron scattering experiment would 
measure the proportion of neutrons that scattered off with the energies 
$m_1, \ldots, m_8$: the resulting $\mathcal{S}(0,\omega)$ would look 
as in the schematic Fig. \ref{fig:mass.peaks}.
Metaphorically speaking, the crystal would thus ``sound'' 
as a ``chord'' of eight ``notes'': the eigenfrequencies $m_i$.

At the lowest temperatures, and in the immediate vicinity of $h_x = h_c$, 
the experiment \cite{Coldea} succeeded to resolve the first two excitations, 
and to extract their masses $m_1$ and $m_2$. 
The mass ratio $m_2/m_1$ was found to be $\frac{m_2}{m_1} = 1.6 \pm 0.025$, 
consistent with $\frac{m_2}{m_1} = \frac{1 + \sqrt{5}}{2} \approx 1.618$ of the 
expression for the $v_{Zam}(m)$ in the Subsection \ref{sub492}. 
In other words, the experimentalists were able to hear two of the 
eight notes of the Zamolodchikov $E_8$ chord. 

A reader wishing to find out more about various facets of the story is 
invited to turn to the references \cite{Rajaraman,Delfino,Goss,Borthwick}.

\begin{acknowledgements}
We are grateful to Misha Finkelberg, Andrei Gabrielov, and Sabir Gusein-Zade 
for the inspiring correspondence, and to Patrick Dorey for sending us his thesis. 
Our special gratitude goes to Sasha Givental whose remarks enabled us to 
generalize some statements and to simplify the exposition. A.V.  
thanks MPI in Bonn for hospitality; he was supported in part by NSF grant 
DMS-1362924 and the Simons Foundation grant  no. 336826.

\end{acknowledgements}

\end{document}